\documentclass[reqno,english]{amsart}
\usepackage{amsfonts,amsmath,latexsym,verbatim,amscd,mathrsfs,color,array}

\usepackage{amsmath,amssymb,amsthm,amsfonts,graphicx,color}
\usepackage{amssymb}
\usepackage{pdfsync}
\usepackage{epstopdf}

\makeatletter
\def\theequation{\thesection.\@arabic\c@equation}
\makeatother

\newcommand{\inn}{{\quad\hbox{in } }}
\newcommand{\ass}{{\quad\hbox{as } }}
\newcommand{\onn}{{\quad\hbox{on } }}
\newcommand{\ttt}{\tilde }

\newcommand{\nn}{ {\nabla}  }

\newcommand{\pp}{ {\partial} }

\newcommand{\A}{\alpha }
\newcommand{\vp}{\varphi}

\newcommand{\R} {\mathbb R}
\newcommand{\La} {\Lambda }
\newcommand{\cuad}{{\sqcap\kern-.68em\sqcup}}

\newcommand{\ppp }{{q}}

\newcommand{\foral}{\quad\mbox{for all}\quad}
\newcommand{\ve}{\varepsilon}

\newcommand{\be}{\begin{equation}}
\newcommand{\ee}{\end{equation}}

\newcommand{\la}{\lambda}

\newcommand{\equ}[1]{(\ref{#1})}

\renewcommand{\theequation}{\thesection.\arabic{equation}}
 
 \newtheorem{lemma}{Lemma}[section]

\newtheorem{teo}{Theorem}

\newtheorem{prop}{Proposition}[section]
\newtheorem{corollary}{Corollary}[section]
\newtheorem{remark}{Remark}[section]
\newcommand{\bremark}{\begin{remark} \em}
\newcommand{\eremark}{\end{remark} }

\begin{document}

\title[Standing waves for the fractional laplacian]
{Concentrating standing waves for the fractional nonlinear Schr\"odinger equation}


\author[J. D\'avila]{Juan D\'avila}
\address{\noindent   Departamento de
Ingenier\'{\i}a  Matem\'atica and Centro de Modelamiento
 Matem\'atico (UMI 2807 CNRS), Universidad de Chile,
Casilla 170 Correo 3, Santiago,
Chile.}
\email{jdavila@dim.uchile.cl}

\author[M. del Pino]{Manuel del Pino}
\address{\noindent   Departamento de
Ingenier\'{\i}a  Matem\'atica and Centro de Modelamiento
 Matem\'atico (UMI 2807 CNRS), Universidad de Chile,
Casilla 170 Correo 3, Santiago,
Chile.}
\email{delpino@dim.uchile.cl}

\author[J. Wei]{Juncheng Wei}
\address{\noindent  Department of Mathematics, Chinese
University of Hong Kong, Shatin, Hong Kong and Department of Mathematics, University of British Columbia, Vancouver, B.C., Canada, V6T 1Z2.
} \email{wei@math.cuhk.edu.hk}

\begin{abstract}
We consider
the semilinear equation
$$
\ve^{2s} (-\Delta)^s u + V(x)u - u^p = 0, \quad u>0, \quad u\in H^{2s}(\R^N)
$$
where $0<s<1,\  1<p<\frac{N+2s}{N-2s}$,
$  V(x)$ is a sufficiently smooth potential with
$\inf_\R V(x) > 0$, and $\ve>0$ is a small number. Letting $w_\la$ be the radial ground state of $(-\Delta )^sw_\la  + \la w_\la - w_\la^p=0$ in $H^{2s}(\R^N)$, we build solutions
of the form  $$ u_\ve(x) \sim \sum_{i=1}^k w_{\lambda_i}  ( (x-\xi_i^\ve)/\ve),$$
where  $\lambda_i = V(\xi_i^\ve)$ and the $\xi_i^\ve $ approach suitable critical points of $V$. Via a Lyapunov Schmidt variational reduction, we recover various existence results already known for the case $s=1$. In particular such a solution exists around $k$ nondegenerate critical points of $V$. For $s=1$ this corresponds to the classical results  by Floer-Weinstein \cite{fw} and Oh \cite{oh1, oh2}.

\end{abstract}

\date{}\maketitle


\color{black}

\setcounter{equation}{0}
\section{Introduction and main results}
We consider the fractional nonlinear Schr\"odinger equation
\be
i\psi_t =   \ve^{2s} (-\Delta)^s\psi + W(x) \psi - |\psi|^{p-1}\psi
\label{nls}\ee
where $(-\Delta)^s$, $0<s<1$, denotes the usual fractional Laplace operator, $W(x)$ is a bounded potential, and $p> 1$.  We are interested in the {\em semi-classical limit} regime,
$0<\ve\ll 1$.

We want to find {\em standing-wave solutions}, which are those of the form $\psi(x,t) = u(x) e^{iE t}$ with $u$  real-valued function.
Letting $V(x) = W(x) + E $, equation \equ{nls} becomes
\be
\ve^{2s} (-\Delta)^su + V(x)u - |u|^{p-1}u = 0 \inn \R^N.
 \label{fw1}
 \ee
We assume in what follows that $V$ satisfies
\be
 V \in C^{1,\alpha}(\R^N)\cap L^\infty (\R^N), \quad  \inf_{\R^N} V(x)  >0.
\label{VV}\ee
We are interested in finding solutions with a {\em spike pattern} concentrating around a finite number of points in space as $\ve\to 0$ .
This has been the topic of many works in the standard case $s=1$, relating the concentration points with critical points of
 the potential, starting in 1986 with the pioneering work by Floer and Weinstein \cite{fw}, then continued by Oh \cite{oh1, oh2}.
 The natural place to look for solutions to \equ{fw1} that decay at infinity is the space $H^{2s}(\R^N)$, of all functions $u\in L^2(\R^N)$ such that
 $$
 \int_{\R^N} (1+ |\xi|^{4s})\, |\hat u (\xi)|^2 d\xi < +\infty,
 $$
 where $\widehat {\, }$ denotes Fourier transform. The fractional Laplacian $(-\Delta)^su$ of a function  $u\in H^{2s}(\R^N)$
  is defined in terms of its Fourier transform by the relation
 $$
 \widehat{(-\Delta)^su } = |\xi|^{2s} \hat u \in L^2(\R^N).
   $$

   We will explain next what we mean by a {\em spike pattern} solution of equation \equ{fw1}.
Let us consider the basic problem
\be
(-\Delta)^sv  + v - |v|^{p-1}v = 0,\quad v\in H^{2s}( \R^N).
\label{fw2}\ee
We assume the following constraint in $p$,
\be\label{constp}
 1< p<  \left \{ \begin{matrix} \frac{ N+ 2s}{ N-2s}  & \hbox{ if }  2s< N, \\ + \infty   & \hbox{ if }  2s\ge N. \end{matrix} \right .\ee
Under this condition it is known the existence of a positive, radial {\em least energy solution}
$v=w(x)$, which gives the lowest possible value for the energy
   $$
J(v) =
\frac 12 \int_{\R^{N}} v (-\Delta)^s v   + \frac 1 2  \int_{\R^N} v^2 - \frac 1{p+1}\int_{\R^N} |v|^{p+1}.
$$
among all nontrivial solutions of \equ{fw2}. An important property, which has only been proven recently by Frank-Lenzman-Silvestre \cite{fls} (see also \cite{at,fl}), is that there exists a radial least energy solution which is nondegenerate, in the sense that the space of solutions of the equation
\be
(-\Delta )^s \phi +\phi - p w^{p-1} \phi=0, \quad  \phi\in H^{2s}(\R^n)
\ee
consists of exactly  of the linear combinations of the translation-generators, $ \frac{\partial w}{\partial x_j},\  j=1, \cdots, N$.

\bigskip
It is easy to see that the function
$$
w_\la (x):= \la^{\frac 1{p-1}} w ( \la^{\frac 1{2s}}x)
$$
satisfies the equation
$$
(-\Delta)^sw_\la  + \la w_\la - w_\la^p = 0\inn \R^N.
$$
Therefore for any point $\xi\in \R^N$, taking $\la = V(\xi)$, the spike-shape function
\be
u(x)= w_{V(\xi)} \left ( \frac {x-\xi}\ve  \right )
\label{sp}\ee
satisfies
$$
\ve^{2s} (-\Delta)^su+ V(\xi) u -  u^p = 0.  \quad
$$
Since the $\ve$-scaling makes it {\em concentrate} around $\xi$, this function constitutes a good positive approximate solution
to equation \equ{fw1}, namely of

\be
\ve^{2s} (-\Delta)^su + V(x)u -  u^p = 0,\quad
 \label{fw3}
 \ee
 $$
 u>0,\quad u\in H^{2s}( \R^N).
$$

We call a {\em $k$-spike pattern solution} of \equ{fw3} one that looks approximately like a superposition of $k$ spikes like \equ{sp}, namely a solution  $u_\ve$ of the form

\be\label{spike}
u_\ve(x) = \sum_{i=1}^k   w_{V(\xi_i^\ve)} \left ( \frac {x-\xi_i^\ve}\ve  \right )\, + o(1)\quad
\ee
for  points $\xi_1^\ve,\ldots, \xi_k^\ve$,
where $o(1)\to 0$ in $H^{2s}(\R^N)$ as $\ve \to 0$.

\medskip
In what follows  we assume that
$p$ satisfies condition \equ{constp} and $V$ condition \equ{VV}.

\medskip

Our first result concerns the existence of multiple spike solution at separate places in the case of stable critical points.

\begin{teo} \label{teo1}
Let $\Lambda_i \subset \R^N$, $i=1,\ldots, k$, $k\ge 1$   be disjoint bounded open sets in $\R^N$.
Assume that
$$
  \deg(\nn V, \La_i,0) \ne 0 \foral i=1,\ldots, k .
$$
Then for all sufficiently small $\ve$, Problem $\equ{fw3}$ has a solution of the form $\equ{spike}$ where
$\xi_i^\ve \in \Lambda_i$ and
$$ \nn V (\xi_i^\ve) \to 0 \ass \ve \to 0. $$

\end{teo}


An immediate consequence of Theorem \ref{teo1} is the following.
\begin{corollary}\label{coro1}
Assume that $V$ is of class $C^2$. Let $\xi_1^0,\ldots, \xi_k^0$ be  $k$ non-degenerate critical points of $V$, namely
$$
\nn V(\xi_i^0) = 0, \quad D^2V(\xi_i^0) \hbox{ is invertible for all } i=1,\ldots k. 
$$
Then, a $k$-spike solution of $\equ{fw3}$ of the form $\equ{spike}$ with $\xi_i^\ve\to \xi_i^0$
exists.
\end{corollary}

When $s=1$, the result of Corollary \ref{coro1} is due to  Floer and Weinstein \cite{fw} for $N=1$ and $k=1$  and to Oh \cite{oh1, oh2} when $N\ge1$, $k\ge 1$.
 Theorem \ref{teo1} for $s=1$ was proven by Yanyan Li \cite{yyli}.

\begin{remark}\label{remark}  {\em
As the proof will yield, Theorem 1 for $0<s<1$ holds true under the following, more general condition introduced in \cite{yyli}. Let
$\La= \La_1\times\cdots \times \La_k$ and assume that the function
\be \vp( \xi_1, \ldots , \xi_k) = \sum_{i=1}^k V(\xi_i)^\theta, \quad \theta = \frac {p+1}{p-1} - \frac N{2s} >0
\label{vp}
\ee has a {\em stable critical point situation} in $\La$: there is a number $\delta_0>0$ such that for each $g \in C^1(\bar \Omega)$ with
$\|g\|_{L^\infty(\Lambda)} + \|\nn g\|_{L^\infty(\Lambda)}< \delta_0 $,  there is a $\xi_g\in\Lambda$ such that $ \nn\vp(\xi_g) + \nn g (\xi_g) = 0 $.
Then for all sufficiently small $\ve$, Problem $\equ{fw3}$ has a solution of the form $\equ{spike}$ where
$\xi^\ve= (\xi_1^\ve,\ldots, \xi_k^\ve) \in \La$ and $\nn \vp(\xi^\ve) \to 0$ as $\ve \to 0$.}
\end{remark}

\bigskip

\begin{teo}\label{teo2}
Let  $\La$ be a bounded, open set with smooth boundary such that $V$ is such that either
\be
c= \inf_\Lambda  V < \inf_{\pp\Lambda}  V
\label{mini}\ee
or
$$
c= \sup_\Lambda  V > \sup_{\pp\Lambda}  V
$$
or,
there exist closed sets $B_0\subset B \subset \La$ such that
\be
c=  \inf_{\Phi\in \Gamma}  \sup_{x\in B} V(\Phi(x))  > \sup_{B_0} V
\label{mini3}. \ee
where $\Gamma =\{ \Phi\in C(B, \bar \La) \ /\  \Phi\Big|_{B_0} = Id\ \} $ and
$\nn V(x)\cdot \tau \ne 0 $ for all $x\in \pp\La$ with $V(x)=c$ and some tangent vector $\tau$ to $\pp\La$ at $x$.

\medskip
Then, there exists a $1$-spike solution of $\equ{fw3}$ with $\xi^\ve\in \La$ with\\ $\nn V(\xi_\ve)\to 0$ and $V(\xi^\ve_i)\to c$.
\end{teo}

In the case $s=1$, the above results  were found by del Pino and Felmer \cite{df1,df2}. The case of a (possibly degenerate) global minimizer was previously considered by Rabinowitz \cite{rabinowitz} and X. Wang \cite{wang}.
 An isolated maximum with a power type degeneracy appears in Ambrosetti, Badiale and Cingolani \cite{abc}. Condition \equ{mini3} is called a {\em nontrivial linking situation} for $V$.
The cases of $k$ disjoint sets where \equ{mini} holds was treated in \cite{df3,gui}. Multiple spikes for disjoint nontrivial linking regions were first considered in \cite{df4}, see also \cite{cingolanilazzo,grossi} for other multiplicity results.

\medskip

Our last result concerns the existence of multiple spikes {\em at the same point}.

\begin{teo}\label{teo3}
Let  $\La$ be a bounded, open set with smooth boundary such that $V$ is such that
$$
 \sup_\Lambda  V > \sup_{\pp\Lambda}  V.
$$

Then for any positive integer $k$ there exists a $k$-spike solution of $\equ{fw3}$ with  spikes $\xi_j^\ve\in \La$ satisfying $ V(\xi_j^\ve)\to \max_{\Lambda} V$.

\end{teo}

In the case $s=1$, Theorem \ref{teo3} was proved by Kang and Wei \cite{kw}.  D'Aprile and Ruiz \cite{daprile} have found a phenomenon of this type
at a saddle point of $V$.

\bigskip
The rest of this paper will be devoted to the proofs of Theorems 1--3. The method of construction of a $k$-spike solution
consists of a Lyapunov-Schmidt reduction in which the full problem is reduced to that of finding a critical point $\xi^\ve$
of a functional which is a small $C^1$-perturbation of  $\vp$ in \equ{vp}. In this reduction
the nondegeneracy result in \cite{fls} is a key ingredient.

\medskip
After this has been done,
 the results follow directly from standard degree theoretical or variational arguments.   The Lyapunov-Schmidt reduction is a method widely
used in elliptic singular perturbation problems. Some results of variational type for $0<s<1$ have been obtained for instance in \cite{fqt} and
\cite{secchi}. We believe that the scheme of this paper  may be generalized to concentration on higher dimensional regions, while that could be much more challenging. See \cite{dkw,mahmoudi} for concentration along a curve in the plane and $s=1$.

\setcounter{equation}{0}
\section{Generalities} Let $0<s< 1$.
Various definitions of the fractional Laplacian $(-\Delta)^s\phi$ of a function $\phi$ defined in $\R^N$ are available, depending on its regularity and growth properties.

\medskip
As we have recalled in the introduction, for $\phi\in H^{2s}(\R^N)$ the standard definition  is given via Fourier transform $\widehat{\ }$.
 $(-\Delta)^s\phi \in L^2(\R^N)$ is defined by the formula
\be\label{s11}
|\xi|^{2s} \hat \phi(\xi) = \widehat { (-\Delta)^s\phi} .
\ee
When $\phi$ is assumed in addition sufficiently regular, we obtain the direct representation
\be
(-\Delta)^s\phi(x) =   \,  d_{s,N} \,\int_{\R^N}  \frac {\phi(x) -\phi(y)}{ |x-y|^{ N+2s} }\, dy
\label{s2}\ee
for a suitable constant $d_{s,N}$ and the integral is understood in a principal value sense.
This integral makes sense directly when $s< \frac 12$ and $\phi \in C^{0,\alpha}(\R^N)$ with $\alpha > 2s$, or if
 $\phi \in C^{1,\alpha} (\R^N) $, $1+\alpha >  2s$. In the latter case,  we can desingularize the integral representing it in the form
$$
(-\Delta)^s\phi(x) =   d_{s,N}\int_{\R^N}  \frac {\phi(x) -\phi(y) -\nn \phi (x)(x-y)  }{ |x-y|^{ N+2s} }\, dy.
$$
Another useful (local) representation, found by Caffarelli and Silvestre \cite{cs}, is via the following
boundary value problem in the half space $ \ \R^{N+1}_+=\{(x,y)\ /\ x\in \R^N,\  y>0\} $:
\begin{equation*}
\left\{\begin{array}{l}
 \nn\cdot  ( y^{1-2s} \nn \ttt \phi  )  =  0  \inn \R^{N+1}_+ , \\
 \quad   \ttt\phi (x,0) = \phi(x)\ \  \onn \R^N.
\end{array}
\right.
\end{equation*}
Here $\ttt\phi$ is the $s$-harmonic extension of $\phi$, explicitly given as a convolution integral with the $s$-Poisson kernel $p_s(x,y)$,
$$
\ttt \phi (x, y)  =  \int_{\R^N}  p_s(x-z,y ) \phi(z) \, dz,
$$
where
$$
\quad p_s (x,y) = c_{N,s}  \frac {y^{4s-1}}
{(|x|^2 +  |y|^2 ) ^{\frac {N-1 + 4s } 2}}
$$
and $c_{N,s}$ achieves $\int_{\R^N} p(x,y)dx =1.$
Then under suitable regularity, $(-\Delta)^s \phi$ is the Dirichlet-to-Neumann map for this problem, namely
\be
(-\Delta)^s \phi (x)  \ =\ \lim_{y\to 0^+} y^{1-2s} \pp_y \ttt\phi (x,y).
\label{s3}\ee
Characterizations \equ{s11}, \equ{s2}, \equ{s3} are all equivalent for instance
in Schwartz's space of rapidly decreasing smooth functions.

\medskip
Let us consider now for a number $m>0$ and $g\in L^2(\R^N)$ the equation
$$
(-\Delta)^s \phi  + m\phi  = g \quad \hbox{in } \R^N.
$$
Then in terms of Fourier transform,  this problem, for $\phi \in  L^2$, reads
$$
(|\xi|^{2s} + m)\, \hat \phi    =  \hat g
$$
and has a unique solution $\phi \in H^{2s}(\R^N)$ given by the convolution
\be
\phi(x) =  T_m[g]:=  \int _{\R^N} k(x-z)\, g(z)\, dz,
\label{poto}\ee
where
$$
\hat k(\xi) = \frac 1{ |\xi|^{2s} + m }.
$$
Using the characterization \equ{s3} written in weak form, $\phi$ can then be characterized by $\phi(x) = \ttt\phi(x,0)$ in trace sense,
where   $\ttt \phi\in H $ is the unique solution of
\be
\iint_{\R^{N+1}_+} \nn \ttt \phi \nn \vp \,y^{1-2s}  +  m \int_{\R^N} \phi\vp  = \int_{\R^N} g\vp, \quad \foral \vp\in H,
\label{pico2}\ee
 where $H$ is the Hilbert space of functions $\vp \in H^1_{loc} (\R^{N+1}_+)$ such that
$$
\|\vp\|_H^2  :=  \iint_{\R^{N+1}_+} |\nn \vp|^2 y^{1-2s} +  m \int_{\R^N} |\vp|^2 \ <\ +\infty,
$$
or equivalently the closure of the set of all functions in $C_c^\infty ( \overline{\R^{N+1}_+}) $ under this norm.

A useful fact for our purposes is the equivalence of the representations \equ{poto} and \equ{pico2} for $g\in L^2(\R^N)$.

\begin{lemma}
Let $g\in L^2(\R^N)$. Then the unique solution  $\ttt \phi\in H$ of Problem \equ{pico2} is given by the $s$-harmonic extension
of the function $\phi = T_m[g] = k* g $.

\end{lemma}

\proof Let us assume first that $\hat g \in C_c^\infty (\R^N)$. Then  $\phi$ given by \eqref{poto} belongs to $H^{2s}(\R^N)$. Take a test function $\psi\in  C_c^\infty (\R^{N+1}_+)$. Then the well-known computation by Caffarelli and Silvestre shows that
$$
\iint_{\R^{N+1}_+} \nn \ttt\phi \nn\psi \, y^{1-2s}dydx  =
$$
$$
\int_{\R^N} \lim_{y\to 0 } y^{1-2s}\pp_y \ttt\phi (y, \cdot)\,\psi\, dx = \int_{\R^N} \psi (-\Delta )^s \phi\, dx = \int_{\R^N} (g-  m\phi) \, dx.
$$
By taking $\psi = \ttt\phi \eta_R$ for a suitable sequence of smooth cut-off functions equal to one on expanding balls $B_R(0)$ in $\R^{N+1}_+$, and using the behavior at infinity of  $\ttt\phi $ which resembles the Poisson kernel $p_s(x,y),$ we obtain
$$
\iint_{\R^{N+1}_+}  |\nn \ttt\phi|^2 \, y^{1-2s}dydx  + m \int_{\R^N} |\phi|^2 =   \int_{\R^N} g\phi
$$
and hence $\|\ttt\phi\|_H \le C\|g\|_{L^2}$ and satisfies \equ{pico2}. By density, this fact extends to all $g\in L^2(\R^N)$.
The result follows since the solution of Problem \equ{pico2} in $H$ is unique. \qed

\medskip
 Let us recall the main properties of the fundamental solution $k(x)$ in the representation \equ{poto}, which are stated for instance in \cite{fls} or
in \cite{fqt}.

We have that  $k$ is radially symmetric and positive,  $k\in C^\infty (\R^N\setminus\{0\})$ satisfying

\begin{itemize}

\item $$ |k(x)| + |x|\,|\nn k(x)| \ \le \ \frac C{ |x|^{N-2s}} \foral |x| \le 1, $$

\item   $$ \lim_{|x|\to \infty }  k(x) |x|^{N+2s}  = \gamma > 0, $$

\item $$|x|\, |\nn k(x) |   \ \le \ \frac C{ |x|^{N+2s}} \foral |x| \ge 1.$$

\end{itemize}

The operator $T_m$ is not just defined on functions in $L^2$. For instance it acts nicely on bounded functions. The positive kernel $k$ satisfies $\int_{\R^N} k =\frac 1m $. We see that if $g\in L^\infty(\R^N)$ then
$$
\|T_m[g]\|_\infty \le \frac 1m \|g\|_\infty.
$$

We have indeed the validity of an estimate like this for $L^\infty$ weighted norms as follows.

\begin{lemma} \label{pico10}
Let $0\le \mu < N +2s$. Then there exists a $C>0$ such that
$$
 \|(1+ |x|)^{\mu}T_m[g]\|_{L^\infty(\R^N)} \le    C\|(1+ |x|)^{\mu}g\|_{L^\infty(\R^N)}.
$$
\end{lemma}

\proof
Let us assume that $0\le \mu < N +2s$ and let $\bar g( x) = \frac 1{ (1+ |x|)^\mu}$.
Then
$$
 T[\bar g] (x) = \int_{|y-x| < \frac 12 |x| } \frac {  k(y) }{  (1+ |y-x|)^\mu}\, dy  + \int_{|y-x| > \frac 12|x| }  \frac {k(y) }{  (1+ |y-x|)^\mu}\, dy.
$$
Then, as $|x|\to \infty $ we find
$$
|x|^\mu \int_{|y-x| < \frac 12 |x| } \frac {k(y) }{  (1+ |y-x|)^\mu}\, dy \sim |x|^{-2s} \to 0,
$$
and since $k\in L^1(\R^N)$, by dominated convergence we find that as $|x|\to \infty $
$$
\int_{|x-y| > \frac 12|x| }  \frac {k(y)|x|^\mu }{  (1+ |x-y|)^\mu}\, dy \to   \int_{\R^N}  k(z)dz = \frac 1m.
$$
We conclude in particular that for a suitable constant $C>0$,  we have
$$
T_m[ (1+ |x|)^{-\mu} ]  \le C(1+ |x|)^{-\mu}.
$$
Now, we have that
$$
\pm T_m[g] \le   \|(1+ |x|)^{\mu}g\|_{L^\infty(\R^N)}\,  T_m[ (1+ |x|)^{-\mu} ],
$$
and then

$$
 \|(1+ |x|)^{\mu}T[g]\|_{L^\infty(\R^N)} \le    C\|(1+ |x|)^{\mu}g\|_{L^\infty(\R^N)}
$$
as desired. \qed

\medskip
We also have the validity of the following useful  estimate.
\begin{lemma}
Assume that $g\in L^2\cap L^\infty$.
Then the following holds: if $\phi = T_m [ g]$ then there is a $C>0$ such that
\be
\sup_{x\ne y} \frac{|\phi(x) -\phi(y) |}{|x-y|^\alpha}  \le  C \|g\|_{L^\infty(\R^N)}
\label{holder}\ee
where $\alpha = \min\{1,2s\}$.
\end{lemma}
\proof
Since  $\|T_m[g]\|_\infty \le C\|g\|_\infty$, it suffices to establish \equ{holder} for $|x-y|< \frac 13$.
We have
$$
|\phi(x) -\phi(y) | \le \int_{\R^N} |k(z+ y-x) - k(z)|\, dz\,  \|g\|_\infty.
$$
Now, we decompose
$$
\int_{\R^N} |k(z+ y-x) - k(z)|\, dz\, =
$$
$$
\int_{|z| > 3|y-x|} |k(z + y-x )-k(z)|\, dz  +  \int_{|z| < 3|y-x|} |k(z + y-x )-k(z)|\, dz.
$$
We have
$$
 \int_{|z| > 3|y-x|} | k(z + (y-x) )-k(z)| \le    \int_0^1dt \int_{|z| > 3|y-x|}  |\nn k( z + t(y-x))|\, dz \, |y-x| .
$$
and, since $3|y-x| < 1$,
$$
 \int_{|z| > 3|y-x|}  |\nn k( z + t(y-x))| dz \le   C( 1 + \int_{1> |z| > 3|y-x|}  \frac {dz} {|z|^{N+1 - 2s}} )\, \le  C(1+ |y-x|^{2s -1}).
$$
On the other hand
$$
 \int_{|z| < 3|y-x|} |k(z + y-x )-k(z)|\, dz \le   2\int_{|z | < 4|y-x|}  |k(z)|\, dz \le  C|y-x|^{2s},
$$
and \equ{holder} readily follows.\qed

\medskip

\bigskip
Next we consider the more general problem
\be
(-\Delta)^s \phi  + W(x)\phi  = g \quad \hbox{in } \R^N
\label{pico3}\ee
where $W$ is a bounded potential.

We start with a form of the weak maximum principle.

\begin{lemma}\label{potin1}
Let us assume that
$$
\inf_{x\in \R^N }  W(x) =: m  > 0
$$
and that $\phi \in H^{2s}(\R^N)
$ satisfies equation $\equ{pico3}$ with $g\ge 0$. Then $\phi \ge 0$ in $\R^N$.
\end{lemma}
\proof
We use the representation for $\phi$ as the trace of the unique solution $\ttt\phi\in H$  to the problem
$$
\iint_{\R^{N+1}_+} \nn \ttt \phi \nn \vp y^{1-2s}  +   \int_{\R^N}  W \phi\vp  = \int_{\R^N}  g\vp, \quad \foral \vp\in H.
$$
It is easy to check that the test function $\vp = \phi_-= \min\{\phi, 0\}$ does indeed belong to $H$. We readily obtain
$$
\iint_{\R^{N+1}_+} |\nn \ttt \phi_-|^2 y^{1-2s}  +   \int_{\R^N}  W\phi_-^2  = \int_{\R^N} g\phi_- .
$$

Since  $g\ge 0$ and $ W\ge m$, we obtain that $\phi_-\equiv 0$, which means precisely $\phi \ge 0$, as desired. \qed

\medskip
We want to obtain a priori estimates for problems of the type \equ{pico3} when $W$ is not necessarily positive.
Let $\mu > \frac N2$, and let us assume that
$$\|(1+ |x|^\mu) g\|_{L^\infty(\R^N)} < +\infty.$$
The assumption in $\mu$ implies that $g\in L^2(\R^N)$.

\medskip
Below, and in all what follows, we will say that $\phi\in L^2(\R^N)$ solves equation \equ{pico3} if and only if $\phi$ solves
the linear problem
$$
\phi =  T_m( (m-W)\phi + g).
$$
Similarly, we will say that
$$
(-\Delta)^s \phi  + W(x)\phi  \ge  g \quad \hbox{in } \R^N
$$
if for some $\ttt g\in L^2(\R^N)$ with $\ttt g \ge g$ we have
$$
\phi =  T_m( (m-W)\phi + \ttt g).
$$

\medskip
The next lemma provides an a priori estimate
for a solution $\phi \in L^2(\R^N)\cap L^\infty (\R^N)$ of \equ{pico3}.

\begin{lemma}\label{lemin}
Let $W$ be a continuous function, such that for $k$ points $q_i$ $i=1,\ldots ,k$  a number $R>0$
and
$
B = \cup_{i=1}^k B_R(q_i)
$
we have
$$
\inf_{x\in \R^N \setminus B}  W(x) =: m  > 0.
$$

Then, given any number $ \frac N2<   \mu < N + 2s$ there exists a constant $C= C(\mu,k,R)>0$ such that for any  $\phi\in  H^{2s}\cap L^\infty (\R^N)$ and
$g$ with  $$\|\rho^{-1} g\|_{L^\infty(\R^N)} < +\infty$$ that satisfy equation \equ{pico3}
we have the validity of the estimate
$$
\| \rho^{-1}  \phi \|_{L^\infty (\R^N )} \ \le \ C\, \left [\,
\|  \phi \|_{L^\infty (B)}  + \| \rho^{-1} g \|_{L^\infty (\R^N)}
 \,\right ].
$$
Here $$\rho(x) =   \sum_{i=1}^k\frac 1{ (1+ |x-q_i|) ^\mu}.$$
\end{lemma}

\proof
We start by noticing that $\phi$ satisfies the equation
$$
(-\Delta )^s\phi   +  \hat W\phi =  \hat g
$$
where
$$\hat g  =   (m- W) \chi_{B}\, \phi,\quad   \hat W =   m\chi_{B}  + W(1-\chi_{B}). $$
Observe that
$$
|\hat g(x) | \le M \sum_{i=1}^k (1+ |x-q_i|)^{-\mu},  \quad  M = C( \|  \phi \|_{L^\infty (B)}  + \| \rho^{-1}  g \|_{L^\infty (\R^N)} )$$
where $C$ depends only on $R$, $k$ and $\mu$ and
$$
\inf_{x\in \R^N} \hat W(x) \ge m .
$$
 Now, from Lemma \ref{pico10}, since  $0<\mu< N+ 2s$ we  find a solution $\phi_0 (x)$
to the problem $$(-\Delta)^s\bar \phi  + m  \bar \phi = (1+ |x|)^{-\mu}$$
such that $\bar \phi = O( |x|^{-\mu})$ as $|x|\to \infty $.
Then we have that
$$((-\Delta)^s + \hat W )  (\bar\phi)  \ge M \sum_{i=1}^k (1+ |x-q_i|)^{-\mu}$$
where
$$
\bar\phi(x) = M\sum_{i=1}^k \phi_0( x-q_i).
$$
Setting $\psi = (\phi - \bar\phi ) $ we get
$$(-\Delta)^s\psi  + \hat W  \psi =  \tilde g \le 0$$
with $\ttt g \in L^2$. Using Lemma \ref{potin1} we obtain $\phi \le \bar \phi$. Arguing similarly for $-\phi$, and
using the form of $\bar\phi $ and $M$, the desired estimate immediately follows. \qed

\medskip
Examining the proof above, we obtain immediately the following.
\begin{corollary}\label{orto}
Let $\rho(x)$ be defined as in the previous lemma. Assume that $\phi\in H^{2s}(\R^N)$ satisfies equation \equ{pico3}
and that
$$
\inf_{x\in \R^N }  W(x) =: m  > 0.
$$
Then we have that $\phi\in L^\infty(\R^N)$ and it satisfies
\be
\| \rho^{-1}  \phi \|_{L^\infty (\R^N )} \ \le \ C\,
 \| \rho^{-1} g \|_{L^\infty (\R^N)}
 \, .
\ee

\end{corollary}

\bigskip
A last useful fact is that if $f,g\in L^2(\R^N)$ and
$
W= T(f), \ Z = T(g)
$
then the following holds:
$$
\int_{\R^N} Z (-\Delta)^s W  -\int_{\R^N} W (-\Delta)^sZ  =  \int_{\R^N} T_m[f]g -\int_{\R^N} T_m[g] f  = 0  ,
$$
the latter fact  since the kernel $k$ is radially symmetric.

\setcounter{equation}{0}
\section{Formulation of the problem: the ansatz}

By a solution of the problem
$$
\ve^{2s}(-\Delta)^s u + V(x) u - u^p = 0 \inn \R^N
$$
we mean a $u\in H^{2s}(\R^N)\cap L^\infty(\R^N)$ such that the above equation is satisfied.
Let us observe that it suffices to solve

\be
\ve^{2s}(-\Delta)^s u + V(x) u - u_+^p = 0 \inn \R^N
\label{piquito00}\ee
where $u_+ =\max\{u,0\}$. In fact, if $u$ solves \equ{piquito00} then
$$
\ve^{2s}(-\Delta)^s u + V(x) u \ge  0 \inn \R^N
$$
and, as a consequence to Lemma \ref{potin1}, $u\ge 0$.

After absorbing $\ve$ by scaling, the equation takes the form
\be
(-\Delta)^s v + V(\ve x) v - v_+^p = 0 \inn \R^N
\label{piquito}\ee
Let us consider points $\xi_1,\ldots, \xi_k\in \R^N $ and designate
$$q_i = \ve^{-1}\xi_i, \quad q= (q_1,\ldots, q_k).$$

Given  numbers   $\delta >0$  small and $R>0$ large, we define the configuration space $\Gamma$ for the points $q_i$ as
\be\label{const}
 \Gamma := \{ q= (q_1,\ldots, q_k)\ /\  R \le \max_{i\ne j} |\ppp_i- \ppp_j|, \quad \max_i |\ppp_i| \le  \delta^{-1}\ve^{-1}\}.
\ee
We look for a solution with concentration behavior near each $\xi_j$. Letting $\ttt v(x) = v(x+ \xi_j)$ translating the origin to  $q_j$,  Equation \equ{piquito} reads
$$
(-\Delta)^s \ttt v + V(\xi_j+ \ve x) \ttt v - \ttt v^p_+ = 0 \inn \R^N.
$$
Letting formally $\ve\to 0$ we are left with the equation
$$
(-\Delta)^s \ttt v + \la_j \ttt v - \ttt v^p_+ = 0 \inn \R^N, \quad \la_j= V(\xi_j).
$$
So we ask that  $v(x) \approx w_{\la_j} (x-q_j ) $ near $q_j$.  We consider the sum of these functions as a first approximation. Thus, we look for a solution $v$ of \equ{piquito}
of the form
$$
v = W_q + \phi, \quad W_q(x) =  \sum_{j=1}^k w_j(x),  \quad w_j(x)= w_{\la_j} (x- q_j) , \quad \la_j = V(\xi_j),
$$
where $\phi$ is a small function, disappearing as $\ve\to 0$.
In terms of $\phi$, Equation \equ{piquito} becomes
\be
(-\Delta)^s \phi + V(\ve x) \phi  - pW_q^{p-1}\phi  = E + N(\phi)       \inn \R^N
\label{piquiton}\ee
where
$$
N(\phi):= (W_q +\phi)_+^p - pW_q^{p-1}\phi - W_q^p  , \quad
$$
\be
 E:= \sum_{j=1}^k (\la_j - V(\ve x)) w_j   + \Big ( \sum_{j=1}^k w_j \Big )^p -  \sum_{j=1}^k w_j^p .
\label{piq3} \ee

Rather than solving Problem \equ{piquiton} directly, we consider first a projected version of it.
Let us consider the functions
$$
Z_{ij}(x) :=  \pp_j w_i(x)
$$
and the problem of finding $\phi \in H^{2s}(\R^N)\cap L^\infty(\R^N)$ such that for certain constants $c_{ij}$

\be
(-\Delta)^s \phi + V(\ve x) \phi  - pW_q^{p-1}\phi  = E + N(\phi) +  \sum_{i=1}^k \sum_{j=1}^N    c_{ij} Z_{ij},
\label{pr1}\ee
\be
   \quad \int_{\R^N} \phi Z_{ij} = 0 \foral i,j.
\label{pr2}\ee

Let ${\mathcal Z } $ be the linear space spanned by
the functions $Z_{ij}$, so  that equation \equ{pr1} is equivalent to
$$
  (-\Delta)^s \phi + V(\ve x) \phi  - pW_q^{p-1}\phi  - E - N(\phi) \in {\mathcal Z }.
$$
On the other hand, for all $\ve$ sufficiently small, the functions $Z_{ij}$ are linearly independent, hence the constants $c_{ij}$ have unique, computable expressions in terms of $\phi$.
We will prove that Problem \equ{pr1}-\equ{pr2} has a unique small solution $\phi= \Phi(q)$. In that way we will get a solution to the full problem \equ{piquiton}
if we can find a value of $q$ such that $c_{ij}(\Phi(q) )=0$ for all $i,j$.
In order to build $\Phi(q)$ we need a theory of solvability for associated  linear operator in  suitable spaces.  This is what we develop in the next section.

\setcounter{equation}{0}
\section{Linear theory}

We consider the linear problem of finding $\phi \in H^{2s} (\R^N) $ such that for certain constants $c_{ij}$ we have
\be
(-\Delta)^s \phi  + V(\ve x) \phi -  pW_\ppp ^{p-1}(x)\phi  +  g(x)  =  \sum_{i=1}^N  \sum_{i=1}^k c_{ij} Z_{ij}
\label{mm1}\ee
\be
\int_{\R^N} \phi Z_{ij} = 0 \foral i,j.
\label{mm2}\ee
The constants $c_{ij} $ are uniquely determined in terms of $\phi$ and $g$ when $\ve$ is sufficiently small, from the linear system

\be
\sum_{i,j} c_{ij} \int_{\R^N}  Z_{ij} Z_{lk}   =   \int_{\R^N} Z_{lk}[ (-\Delta)^s \phi  + V(\ve x) \phi -  pW_\ppp ^{p-1}(x)\phi +g ].
\ee

Taking into account that
$$
\int_{\R^N} Z_{lk} (-\Delta)^s \phi = \int_{\R^N}  \phi (-\Delta)^sZ_{lk}  =  \int_{\R^N} (p w_l^{p-1} -\lambda_l)  Z_{lk} \phi,
$$
  we find
\be
c_{ij} \int_{\R^N}  Z_{ij} Z_{lk}   =   \int_{\R^N}  gZ_{lk} + (p  w_l^{p-1} - pW_\ppp ^{p-1}  + V(\ve x) - \lambda_l) \,  Z_{lk}\phi .
\label{dd}\ee
On the other hand, we check that $$\int_{\R^N}  Z_{ij} Z_{lk} =  \A_l \delta_{ijkl} + O( d^{-N} )   $$ where the numbers $\A_l$ are positive, and independent of $\ve$, and
$$
 d = \min\{ |\ppp _i -\ppp _j |\ /\ i\ne j      \}  \gg 1.
$$
Then, we see that relations \equ{dd} define a uniquely solvable (nearly diagonal) linear system, provided that $\ve$ is sufficiently small.
We assume this last fact in what follows, and hence that the numbers $c_{ij}= c_{ij}(\phi,g)$ are defined by
relations \equ{dd}.

Moreover, we have that
$$
|(p  w_l^{p-1} - pW_\ppp ^{p-1}  + V(\ve x) - \lambda_l) \,  Z_{lk}(x)| \,\le\,  C\,( R^{-N} + \ve |x-q_j|) (1+  |x-q_j|)^{-N-s})
$$
and then from expression \equ{dd} we  obtain the following estimate.

\begin{lemma}\label{estc}
 The numbers $c_{ij}$ in $\equ{mm1}$ satisfy:
$$
c_{ij}  =   \frac 1{\alpha_i} \int_{\R^N} gZ_{ij}    +   \theta_{ij}  .
$$
where
$$
 |\theta_{ij}| \le C (\ve + d^{-N})\, \left [\|\phi\|_{L^2(\R^N)} + \|g\|_{L^2(\R^N)}\right].
 $$
\end{lemma}

In the rest of this section we shall build a solution to Problem \equ{mm1}-\equ{mm2}.

\begin{prop}\label{pro1} Given $k\ge 1$, $\frac N2 < \mu < N+2s$, $C>0$,
there exist positive numbers $d_0$, $\ve_0$, $C$ such that for any points $\ppp _1,\ldots \ppp _k$ and any $\ve$  with
$$ \sum_{i=1}^k|q_i| \le \frac C\ve, \quad  R: = \min\{ |\ppp _i -\ppp _j |\ /\ i\ne j      \} > R_0,  \quad 0< \ve < \ve_0$$
there exists a solution $\phi=T[g]$  of \equ{mm1}-\equ{mm2} that defines a linear operator of $g$, provided that
$$
\|\rho(x)^{-1} g\|_{L^\infty(\R^N) } <+\infty, \quad \rho(x) = \sum_{j=1}^k  \frac 1 {(1 + |x-\ppp _j |)^\mu }.
$$
Besides
$$
\|\rho(x)^{-1} \phi\|_{L^\infty(\R^N) }  \le  C\|\rho(x)^{-1} g\|_{L^\infty(\R^N) }.
$$

\end{prop}

To prove this result we require several steps. We begin with corresponding a priori estimates.

\begin{lemma}  Under the conditions of Proposition \ref{pro1}, there exists a $C>0$ such that
for any solution of \equ{mm1}-\equ{mm2}  with $\|\rho(x)^{-1} \phi\|_{L^\infty(\R^N) } <+\infty$  we have the validity of the a priori estimate
$$
\|\rho(x)^{-1} \phi\|_{L^\infty(\R^N) }  \le  C\|\rho(x)^{-1} g\|_{L^\infty(\R^N) }.
$$

\end{lemma}

\proof
Let us assume the a priori estimate does not hold, namely there are sequences $\ve_n\to 0$, $\ppp _{jn} $, $j=1, \ldots k$,
with $$
 \min\{ |\ppp _{in}-\ppp _{jn} |\ /\ i\ne j      \}\to \infty
$$
and $\phi_n , \ g_n$ with
$$
\|\rho_n(x)^{-1} \phi_n\|_{L^\infty(\R^N) }  =1, \quad \|\rho_n(x)^{-1} g_n\|_{L^\infty(\R^N) } \to 0 ,
$$
where
$$
 \rho_n(x) = \sum_{j=1}^k  \frac 1 {(1 + |x-\ppp _{jn} |)^\mu },
$$
with $\phi_n, g_n$ satisfying \equ{mm1}-\equ{mm2}.
We claim that for any fixed $R>0$ we have that  \be\label{coti} \sum_{j=1}^k \| \phi_n \|_{L^\infty(B_R(\ppp _{jn} ))}\to 0. \ee
Indeed, assume that for a fixed $j$  we have that  $\| \phi_n \|_{L^\infty(B_R(\ppp _{jn} ))}\ge \gamma >0$. Let us set
$ \bar\phi_n (x) = \phi_n ( \ppp _{jn}  + x )$.
We also assume that
 $\la_j^n= V(\ppp _{jn} )\to \bar\la >0$ and
$$
(-\Delta)^s \bar\phi_n    + V(\ppp _{jn} + \ve_n x) \bar \phi_n   + p( w_{\la_j^n}(x) + \theta_n(x))^{p-1}   \bar \phi_n
  = \bar g_n
$$
where
$$
\bar g_n(x) =  g_n ( \ppp _{jn}  + x ) -  \sum_{l=1}^k\sum_{i=1}^n c_{ln}^i\pp_i w_{\la_l^n}(\ppp _{jn} - \ppp _{ln}'+  x ) .
$$
We observe that $\bar g_n(x) \to 0$ uniformly on compact sets.
From the uniform H\"older estimates \equ{holder}, we also obtain equicontinuity of the sequence $\bar\phi_n$. Thus, passing to a subsequence,
we may assume that $\bar\phi_n$ converges, uniformly on compact sets, to a bounded function $\bar \phi$ which satisfies
$\|\bar\phi\|_{L^\infty( B_R(0))} \ge \gamma$.
In addition, we have that
$$
\|(1+|x|)^\mu \bar\phi\|_{L^\infty( \R^N )} \le 1
$$
and that $\bar\phi$ solves the equation
$$
(-\Delta)^s \bar\phi    + \bar \la \bar \phi   + p w_{\bar \la}^{p-1}   \bar \phi
  = 0
  $$
 Let us notice that $\bar\phi \in L^2(\R^N)$, and hence the nondegeneracy result in \cite{fls} applies to yield that
 $\bar\phi$ must be a linear combination of the partial derivatives $\pp_i w_{\bar\la}$.
 But the orthogonality conditions pass to the limit, and yield
 $$
 \int_{\R^N} \pp_i w_{\bar\la}\bar\phi = 0 \foral i=1,\ldots, N.
 $$
 Thus, necessarily $\bar\phi =0$. We have obtained a contradiction that proves the validity of \equ{coti}.
 This and the a priori estimate in Lemma \ref{lemin} shows that also,
$\|\rho_n(x)^{-1} \phi_n\|_{L^\infty(\R^N) }\to 0$, again a contradiction that proves the desired result. \qed

 \medskip
 Next we construct a solution to problem \equ{mm1}-\equ{mm2}.
 To do so, we consider first the auxiliary problem
 \be
(-\Delta)^s \phi +   V\phi  = g +  \sum_{i=1}^k \sum_{j=1}^N    c_{ij} Z_{ij},
\label{prr1}\ee
\be
   \quad \int_{\R^N} \phi Z_{ij} = 0 \foral i,j.
\label{prr2}\ee
where $V$ is our bounded, continuous potential with
$$
\inf_{\R^N} V  = m >0
$$

\begin{lemma}\label{potin}
For each $g$ with $\|\rho^{-1}g\|_\infty <+\infty$, there exists a unique solution of Problem \equ{mm1}-\equ{mm2},
$\phi=: A[g] \in H^{2s}(\R^N)$. This solution satisfies
\be\label{blabla}
\|\rho^{-1}A[g]\|_{L^\infty(\R^N)} \le  C\|\rho^{-1}g\|_{L^\infty(\R^N)}.
\ee
\end{lemma}
\proof
 First we write a variational formulation for this problem.
 Let $X$ be the closed subspace of $H$ defined as
 $$
 X= \{ \ttt \phi\in H \ /\ \int_{\R^N} \phi Z_{ij} = 0 \foral i,j\}
 $$
 Then, given $g\in L^2$, we consider the problem of finding a  $\ttt \phi \in X$ such that
 \be
\langle  \ttt \phi, \ttt\psi   \rangle :=  \iint_{\R^{N+1}_+ } \nn\ttt \phi\nn\ttt \psi y^{1-2s} +  \int_{\R^N} V \phi\psi 
  = \int_{\R^N} g\psi\foral \psi\in X.
 \label{jjj}\ee
We observe that $\langle  \cdot , \cdot   \rangle$
 defines an inner product in $X$ equivalent to that of $H$. Thus existence and uniqueness of a solution follows from Riesz's theorem.
Moreover, we see that
$$ \|\phi\|_{L^2(\R^N) } \le C\|g\|_{L^2(\R^N) }. $$

\medskip
Next we check that this produces a solution in strong sense.
Let $\mathcal Z$ be the space spanned by the functions $Z_{ij}$.  We denote by $\Pi[g]$ the $L^2(\R^N)$ orthogonal projection of $g$ onto
$\mathcal Z$ and by $\ttt \Pi[g]$ its natural $s$-harmonic extension. For a function $\ttt \vp \in H$ let us write
$$
 \ttt \psi = \ttt \vp  - \ttt \Pi[\vp]
$$
so that $\ttt\psi \in X$.
Substituting this $\ttt \psi$ into \equ{jjj} we obtain
$$
 \iint_{\R^{N+1}_+ } \nn\ttt \phi\nn\ttt \vp y^{1-2s} +  \int_{\R^N} V \phi\vp =
$$
$$
\int_{\R^N}  g\vp + \int_{\R^N} [ V \phi  - g ] \, \Pi[\vp]  + \int_{\R^{N}}  \phi (-\Delta )^s\Pi[\vp].
$$
Here we have used that $\ttt\Pi[\vp]$ is regular and
$$
\iint_{\R^{N+1}_+} \nn\phi \nn \ttt \Pi[\vp] y^{1-2s}  =  \int_{\R^{N}}  \phi (-\Delta )^s\Pi[\vp].
$$
Let us observe that for $f\in L^2(\R^N)$ the functional
$$
  \ell (f) = \int_{\R^{N}}  \phi (-\Delta )^s\Pi[f]
$$
satisfies $$|\ell (f)| \le  C\|\phi\|_{L^2(\R^N)} \|\psi\|_{L^2(\R^N)},$$ hence there is an $h(\phi) \in L^2(\R^N)$ such that
$$
\ell (\psi) = \int_{\R^{N}}  h\psi.
$$
If $\phi$ was a priori known to be in $H^{2s}(\R^N)$ we would have precisely that $$h(\phi) = \Pi[ (-\Delta )^s\phi]. $$
Since $\Pi$ is a self-adjoint operator in $L^2(\R^N)$ we then find that
$$
 \iint_{\R^{N+1}_+ } \nn\ttt \phi\nn\ttt \vp y^{1-2s} +  \int_{\R^N} V \phi\vp = \int_{\R^N}  \bar g\vp
$$
where
$$
\bar g = g +  \Pi [ V \phi  - g ] \,   +  h(\phi).
$$
Since $\bar g \in L^2(\R^N)$, it follows then that $\phi \in H^{2s}(\R^N)$ and it satisfies
$$
 (-\Delta )^s\phi + V \phi   - g =   \Pi [ (-\Delta )^s\phi + V \phi  - g ]\in {\mathcal Z},
$$
hence equations \equ{prr1}-\equ{prr2} are satisfied.  To establish estimate \equ{blabla}, we use just Corollary \ref{orto},
observing that
$$
\|\rho^{-1} \Pi [ (-\Delta )^s\phi + V \phi  - g ]\|_{L^\infty(\R^N)} \le  C( \|\phi\|_{L^2(\R^N) } + \|g\|_{L^2(\R^N) })\ \le $$
$$
 C\|g\|_{L^2(\R^N)}\ \le\  \|\rho^{-1} g\|_{L^\infty(\R^N)}.
$$
The proof is concluded. \qed

\medskip

\noindent{\bf Proof of Proposition \ref{pro1}.}
Let us solve now Problem \equ{mm1}-\equ{mm2}. Let $Y$ be the Banach space
\be
Y:= \{ \phi\in C(\R^N) \ /\  \|\phi\|_Y:= \| \rho^{-1}\phi\|_{L^\infty(\R^N)} < +\infty\}
\label{Y}\ee
Let $A$ be the operator defined in Lemma \ref{potin}.
Then we have a solution to Problem \equ{mm1}-\equ{mm2} if we solve
\be
\phi  - A[pW_\ppp ^{p-1} \phi] = A[ g], \quad \phi\in Y.
\label{eee}\ee
We claim that
$$B[\phi] := A[pW_\ppp ^{p-1} \phi] $$
defines  a compact operator in $Y$.
Indeed. Let us assume that $\phi_n$ is a bounded sequence in $Y$. We observe that for some $\sigma>0$ we have
$$
|W_\ppp ^{p-1} \phi_n| \le   C \|\phi_n\|_Y\, \rho^{1+\sigma}.
$$
If $\sigma$ is sufficiently small, it follows that $f_n := B[\phi_n]$ satisfies
$$
|\rho^{-1} f_n| \le  C\rho^\sigma
$$
Besides,
since $f_n = T_m ( (V-m)f_n  +  g_n) $
we use estimate \equ{holder} to get that
for some $\A>0$
$$
\sup_{x\ne y} \frac {|f_n(x)- f_n(y)|}{|x-y|^\A} \le C.
$$
Arzela's theorem then yields the existence of a subsequence of $f_n$
which we label the same way, that converges uniformly on compact sets to a continuous function  $f$ with
$$
|\rho^{-1} f| \le  C\rho^\sigma.
$$
Let $R>0$ be a large number . Then we estimate
$$
\|\rho^{-1}(f_n - f)\|_{L^\infty(\R^N)} \le \|\rho^{-1}(f_n - f)\|_{L^\infty(B_R(0))}  + C\max_{|x|>R}\rho^\sigma(x)
.$$
Since $$\max_{|x|>R}\rho^\sigma(x)\to 0\quad\hbox{ as } R\to \infty$$
 we conclude then that $\|f_n-f\|_\infty \to 0$ and the claim is proven.

\medskip
Finally, the a priori estimate tells us that for $g=0$, equation \equ{eee} has only the trivial solution.
The desired result follows at once from Fredholm's alternative. \qed

\bigskip
We conclude this section by analyzing the differentiability with respect to the parameter $q$ of the solution $\phi= T_q[g]$
of \equ{mm1}-\equ{mm2}. As in the proof above we let $Y$ be the space in \equ{Y}, so that $T_q\in {\mathcal L}(Y)$

\begin{lemma}\label{difer}
The map $q\mapsto T_q$ is continuously differentiable, and for some $C>0$,
\be
\|\pp_q T_q\|_{{\mathcal L}(Y)} \le C
\label{culo}\ee
for all $q$ satisfying constraints \equ{const}.

\end{lemma}

\proof

Let us write $q= (q_1,\ldots, q_k)$, $q_i = (q_{i1},\ldots, q_{iN})$, $\phi = T_q[g]$, and
 (formally) $$\psi = \pp_{q_{ij}}T_q[g], \quad d_{lk} = \pp_{q_{ij}}c_{lk} .$$ Then, by  differentiation of equations \equ{mm1}-\equ{mm2}, we get
 \be
 (-\Delta)^s \psi  + V(\ve x) \psi  - pW_q^{p-1}\psi    =  p\pp_{q_{ij}}W_q^{p-1}\phi +  \sum_{l, k} c_{lk}\pp_{q_{ij}}Z_{lk} +  \sum_{l,k}d_{lk} \,Z_{lk},
\label{fg1}\ee
\be
\int_{\R^N}  \psi Z_{lk}     =    - \int_{\R^N}  \phi \pp_{q_{ij}}Z_{lk}  \foral l,k.
\label{fg2}\ee
We let  $$\ttt \psi = \psi -\Pi[\psi]$$ where, as before, $\Pi[\psi]$ denotes the orthogonal projection of $\psi$ onto the space
spanned by the $Z_{lk}.$
Writing
\be
\Pi[\psi] = \sum_{l,k} \alpha_{lk} Z_{lk}
\label{ffg3}\ee
and relations \equ{fg2} as
\be
\int_{\R^N}  \Pi[\psi ]Z_{lk}     =    - \int_{\R^N}  \phi \pp_{q_{ij}}Z_{lk}  \foral l,k,
\label{ffg2}
\ee
we get
\be
|\alpha_{lk}| \le C \|\phi\|_Y \le C\|g\|_Y.
\label{fg4}
\ee
From \equ{fg1} we have then that
 \be
 (-\Delta)^s \ttt \psi  + V(\ve x) \ttt \psi  - pW_q^{p-1}\ttt \psi    = \ttt g +  \sum_{l,k}d_{lk} \,Z_{lk},
\label{fg3}\ee
or $\ttt \psi= T_q[\ttt g] $
where
\be \ttt g =   p\pp_{q_{ij}}W_q^{p-1}\phi +  \sum_{l,k} c_{lk}\pp_{q_{ij}}Z_{lk} -  [(-\Delta)^s + V(\ve x)   - pW_q^{p-1} ]\,\Pi[\psi] .
\label{fg5}\ee
Then we see that
$$
\|\ttt \psi \|_Y \le C\|\ttt g\|_Y.
$$
Using \equ{fg4} and Lemma \ref{estc}, we see also that
$$
\|\ttt g \|_Y \le C\|g\|_Y, \quad \|\Pi[\psi]\|\le C\|g\|_Y
$$
and thus
\be\label{wq} \|\psi\| \le C\|g\|_Y. \ee
Let us consider now, rigorously, the unique $\psi = \ttt \psi + \Pi[\psi]$ that satisfies equations \equ{fg2} and \equ{fg5}.
We want to show that indeed
$$
\psi = \pp_{q_{ij}}T_q [g].
$$
To do so,
$
q_{i}^t = q_i + te_j
$
where $e_{j}$ is the $j$-th element of the canonical basis of $\R^N$,
and set
$$ q^t =( q_1,\ldots q_{i-1},q_i^t,\ldots, q_k). $$
For a function $f(q)$ we denote
$$
\quad D_{ij}^tf = t^{-1}(f(q^t) -f(q))
$$
we also set
$$
\phi^t:= T_{q^t}[g], \quad   D_{ij}^t T_q[g] =: \psi^t = \ttt\psi^t +  \Pi[ \ttt\psi^t]
$$
so that
$$
 (-\Delta)^s \ttt \psi^t  + V(\ve x) \ttt \psi^t  - pW_q^{p-1}\ttt \psi^t    = \ttt g^t +  \sum_{l,k}d_{lk}^t \,Z_{lk},
$$
where
 $$
 \ttt g^t =   p D_{ij}^{t}[W_{q}^{p-1}]\phi  +  \sum_{l,k} c_{lk}  D_{ij}^t Z_{lk} -  [(-\Delta)^s + V(\ve x)   - pW_q^{p-1} ]\,\Pi[\psi^t] ,$$ $$
\quad d^t_{lk} = D_{ij}^tc_{lk}$$
and
$$
\Pi[\psi^t] = \sum_{l,k} \alpha_{lk}^t Z_{lk} ,
$$
where the constants $\alpha_{lk}^t$ are determined by the relations
$$
\int_{\R^N} \Pi[\psi^t]Z_{lk} = -\int_{\R^N} \phi D_{ij}^{t}Z_{lk} .
$$
Comparing these relations with \equ{ffg3}, \equ{ffg2}, \equ{fg3} defining $\psi$,
we obtain that
$$
\lim_{t\to 0} \|\psi^t - \psi\|_Y = 0
$$
which by definition tells us
$\psi= \pp_{q_{ij}} T_q[g]$. The continuous dependence in $q$ is clear from that of the data in the definition of $\psi$.
Estimate \equ{culo} follows from \equ{wq}. The proof is concluded. \qed

\setcounter{equation}{0}
\section{Solving the nonlinear projected problem}
In this section we solve the nonlinear projected problem
\be
(-\Delta)^s \phi + V(\ve x) \phi  - pW_q^{p-1}\phi  = E + N(\phi) +  \sum_{i=1}^k \sum_{j=1}^N    c_{ij} Z_{ij},
\label{nl11}\ee
\be
   \quad \int_{\R^N} \phi Z_{ij} = 0 \foral i,j.
\label{nl22}\ee

We have the following result.

\begin{prop}\label{prop2}
Assuming that $\|E\|_Y $ is sufficiently small problem $\equ{nl11}$-$\equ{nl22}$
has a unique small solution $\phi = \Phi(q)$ with $$\|\Phi(q)\|_Y \le C\|E\|_Y$$
The map $q\mapsto \Phi(q)$ is of class $C^1$, and for some $C>0$
\be\label{culo1}
\| \pp_q \Phi(q)\|_{ Y} \ \le\ C [ \|E\|_Y + \|\pp_q E\|_Y ]. \ee
for all $q$ satisfying constraints $\equ{const}$.
\end{prop}

\proof
Problem $\equ{nl11}$-$\equ{nl22}$ can be written as
the fixed point problem
\be
\phi = T_q(E + N(\phi)) =:K_q(\phi), \quad \phi \in Y.
\label{fp}\ee
Let
$$
B = \{ \phi \in Y\ /\ \|\phi\|_Y \le  \rho\}..
$$
If $\phi\in B$ we have that
$$
|N(\phi)| \le C|\phi|^\beta, \quad \beta =\min\{p,2\}.
$$
and hence
$$
\|N(\phi)\|_Y \le C\|\phi\|^2
$$
It follows that
$$
\|K_q(\phi) \|_Y \le C_0[ \|E\|   + \rho^2 ]
$$
for a number $C_0$, uniform in $q$ satisfying \equ{const}.
Let us assume
$$
 \rho :=    2C_0 \|E\|,     \quad \|E\| \le  \frac 1{2C_0}   .
 $$
 Then
 $$
 \|K_q(\phi) \|_Y\le C_0[  \frac 1{2C_0} \rho  +   \rho^2 ] \le \rho
 $$
so that $K_q(B) \subset B$. Now, we observe that

$$
|N(\phi_1) - N(\phi_2)| \le  C[|\phi|^{\beta-1} + |\phi|^{\beta-1}]|\phi_1-\phi_2|
$$
and hence
$$
\|N(\phi_1) - N(\phi_2)\|_Y \le C\rho^{\beta-1} \|\phi_1-\phi_2\|_Y
$$
and
$$
\|K_q(\phi_1) - K_q(\phi_2)\| \le C\rho^{\beta-1}\|\phi_1-\phi_2\|_Y.
$$
Reducing $\rho$ if necessary, we obtain that $K_q$ is a contraction mapping and hence has a unique solution of equation \equ{fp} exists in $B$.
We denote it as $\phi = \Phi(q)$.
We prove next that $\Phi$ defines a $C^1$ function of $q$.
Let
$$
M(\phi, q) : = \phi - T_q( E + N(\phi))
$$
Let $\phi_0 = \Phi(q_0)$. Then
$M(\phi_0, q_0)=0$. On the other hand,
$$
\pp_\phi M(\phi, q)[\psi] =  \psi - T_q(N'(\phi)\psi)
$$
where $N'(\phi) = p [( W+ \phi)^{p-1} -W^{p-1}]$,
so that
$$\|N'(\phi)\psi\|_Y \le C\rho^{\beta-1} \|\psi\|_Y $$
If $\rho$ is sufficiently small we have then that $D_\phi M(\phi_0, q_0)$ is an invertible operator, with uniformly bounded inverse.
Besides
$$
\pp_q M(\phi, q) =   (\pp_qT_q)( E + N(\phi)) + T_q( \pp_q E + \pp_q N(\phi))
$$
Both partial derivatives are continuous in their arguments. The implicit function applies in a small neighborhood
of $(\phi_0,q_0)$ to yield existence and uniqueness of a function $\phi= \phi(q)$ with $\phi(q_0)= \phi_0$ defined near $q_0$ with
$
 M(\phi(q), q)=0.
 $
Besides, $\phi(q)$ is of class $C^1$. But, by uniqueness, we must have $\phi(q)=\Phi(q)$.
Finally, we see that
$$
 \pp_q \Phi(q) = -  D_\phi M(\Phi(q), q)^{-1} \left[  (\pp_qT_q)( E + N(\Phi(q))) + T_q( \pp_q E + \pp_q N(\Phi(q)))\right ]
$$
 $$
    \pp_q N(\phi)  =  p[(W+ \phi)^{p-1} - p W^{p-1}   - (p-1)W^{p-2} \phi]\pp_q W   
    $$
and hence
$$
\|(\pp_q N)(\Phi(q))\|_Y \le C\|\Phi(q)\|_Y^\beta \le C \|E\|_Y^\beta
$$
From here,  the above expressions and the bound of Lemma  \ref{difer} we finally get the validity of Estimate \equ{culo1}. \qed

\subsection{An estimate of the error}
Here we provide an estimate of the error $E$ defined in (\ref{piq3}),
 $$ E:= \sum_{j=1}^k (\la_j - V(\ve x)) w_j   + \Big ( \sum_{j=1}^k w_j \Big )^p -  \sum_{j=1}^k w_j^p $$
in the norm $ \| \cdot \|_Y$. Here we need to take $ \mu  \in (\frac{N}{2}, \frac{N+2s}{2})$. We denote
$$ R= \min_{i \not = j} |q_i -q_j |>>1.$$

The first term in $E$ can be easily estimated as
\begin{equation*}
 | \rho^{-1} (x) \sum_{j=1}^k ( \la_j-V(\ve x)) w_j | \leq C \varepsilon^{\min (2s,1) }.
\end{equation*}

To estimate the interaction term in $E$, we divide the $\R^N$ into the $k$ sub-domains
$$ \Omega_j= \{ w_j \geq w_i, \forall i \not = j  \},\quad  j=1, \cdots, k.$$
 In $\Omega_j$, we have
\begin{eqnarray*}
 | (\sum_{j=1}^k w_j  )^p -  \sum_{j=1}^k w_j^p | &\leq & C w_j^{p-1} \sum_{i \not =j} \frac{1}{ |x- q_i|^{N+2s}} \\
& \leq & C \frac{1}{(1+|x-q_j|)^{ (N+2s)(p-1)+ \mu}} \sum_{i \not = j} \frac{1}{ |q_j-q_i|^{N+2s -\mu}}
\\
& \leq & C \rho (x)  R^{ \mu-N-2s}
\end{eqnarray*}

In summary, we conclude that
\begin{equation}
\label{est2}
\| E \|_{Y} \leq C \varepsilon^{2s}+ C R^{\mu-N-2s}
\end{equation}

As a consequence of Proposition \ref{prop2} and the estimate (\ref{est2}), we  obtain that
\begin{equation*}
\|\Phi(q)\|_Y \le C \varepsilon^{\min (2s,1)}+ C R^{\mu-N-2s}.
\end{equation*}

Let us now take
$$ \tau=C \varepsilon^{\min (2s,1) }+ C R^{\mu-N-2s}$$

\setcounter{equation}{0}

\section{The variational reduction}
We will use the above introduced ingredients to find existence results for the equation
\be
(-\Delta)^s v + V(\ve x)v  - v_+^p =0
\label{wz}\ee
An energy whose Euler-Lagrange equation corresponds formally to \equ{wz}
is given by
$$
J_\ve(\ttt v):= \frac 12 \int_{\R^{N}} v(-\Delta)^s v   +  V(\ve x) v^2  -  \frac 1{p+1} \int_{\R^N} V(\ve x) v^2
$$
We want to find a solution of \equ{wz} with the form
$$
v =  v_q:=  W_q  +  \Phi(q)
$$
where $\Phi(q)$ is the function in Proposition \ref{prop2}. We observe that
\be
(-\Delta)^s v_q + V(\ve x) v_q - (v_q)_+^p = \sum_{i,j} c_{ij}Z_{ij}
\label{fr}\ee
hence what we need is to find points $q$ such that $c_{ij}=0$ for all $i,j$. This problem can be formulated variationally
as follows

\begin{lemma}\label{reduction}
Let us consider the function of points $q=(q_1,\ldots, q_k)$ given by
$$
I(q) :=  J_\ve( W_q + \Phi(q)) .
$$
where $ W_q + \Phi(q)$ is the unique $s$-harmonic extension of $W_q +\Phi(q)$.
Then in \equ{fr}, we have $c_{ij}=0$ for all $i,j$ if and only if
$$
\pp_q I(q) = 0.
$$

\end{lemma}

\proof
Let us write $v_q = W_q + \phi(q)$.
We observe that
$$
\pp_{q_{ij}} I(q) =  \int_{\R^{N+1}_+} \nn \ttt v_q \nn(\pp_{q_{ij}}\ttt v_q) y^{1-2s} +  \int_{\R^N} V(\ve x) v_q\pp_{q_{ij}}v_q
-\int_{\R^N} (v_q)_+^{p-1}\pp_{q_{ij}}v_q =
$$
\be
\int_{\R^N}[ (-\Delta)^s v_q + V(\ve x)v_q  - (v_q)_+^p ] \pp_{q_{ij}}v_q = \sum_{k,l} c_{kl} \int_{\R^N} Z_{kl}\pp_{q_{ij}}v_q .
\label{ss}\ee
We observe that
$$
\pp_{q_{ij}}v_q  = -Z_{ij} + O(\ve \rho) +  \pp_{q_{ij}}\Phi(q)
$$
Since, according to Proposition \equ{prop2}
$$ \|\pp_q\Phi(q)\|_Y =  O(\|E\|_Y + \|\pp_q E\|_Y)$$
and this quantity gets smaller as the number $\delta$ in \equ{const} is reduced,
 and the functions $Z_{kl}$ are linearly independent (in fact nearly orthogonal in $L^2$), it follows that the  quantity in \equ{ss} equals zero for all $i,j$ if and only if
$c_{ij}=0$ for all $i,j$. The proof is concluded.  \qed

\bigskip
Our task is therefore to find critical points of the functional $I(q)$. Useful to this end is to achieve
expansions of the energy in special situations.

\medskip

\begin{lemma}
Assume that the numbers $\delta$ and $R$ in the definition of $\Gamma$ in $\equ{const}$
is taken so small that
$$
\|E\|_Y+ \|\pp_q E\| \le  \tau \ll 1.
$$
Then
$$
I_\ve(q)  = J_\ve(W_q) + O(\tau^2)
$$
and
$$
\pp_q I_\ve(q)  = \pp_q J_\ve(W_q) + O(\tau^2)
$$
uniformly on points $q$ in $\Gamma$.

\end{lemma}

\proof
Let us estimate

$$
I(q) =  J_\ve (v_q) ,\quad v_q= W_q + \Phi(q).
$$
 We have that
$$
I(\xi) =   \frac 12 \int_{\R^N} v_q(-\Delta)^s v_q  + V v_q^2  -\frac 1{p+1}\int v_q^{p+1}
 $$
Thus we can expand
$$
I(q) =   J_\ve(W_q)  +   \int_{\R^N} \Phi[(-\Delta)^sv_q + Vv_q - v_q^p]  +
\frac 12 \int_{\R^N} \Phi(-\Delta)^s\Phi  + V\Phi^2
$$
$$
-\frac 1{p+1} \int_{\R^N} [(W_q +\Phi)^{p+1} - W_q^{p+1} - (p+1)W_q^{p}\Phi]
$$
Since, $\|E\|_Y \le \tau $  then  $\|\Phi\|_Y = O(\tau )$, and from the equation satisfied by $\Phi$, also
$\|(-\Delta)^s\Phi\|_Y = O(\tau)$.  This implies
$$|\frac 12 \int_{\R^N} \Phi(-\Delta)^s\Phi  + V\Phi^2| \leq C \int_{\R^N} \rho^{2 \mu}   \tau^2 \leq C \tau^2
$$
and
$$ |\int_{\R^N} [(W_q +\Phi)^{p+1} - W_q^{p+1} - (p+1)W_q^{p}\Phi] | \leq
C \int_{\R^N} \rho^{2 \mu}   \tau^2 \leq C \tau^2
$$

Here we have used the fact that $\mu \in (\frac{N}{2}, \frac{N+2s}{2})$.

On the other hand the second term in the above expansion equals 0,
since by definition $$(-\Delta)^sv_q + Vv_q - v_q^p\in \mathcal{Z}$$ and
$\Phi$ is $L^2$-orthogonal to that space. We arrive to the conclusion that
$$
I(q) = J_\ve (W_q) + O(\tau^2)
$$
uniformly for $q$ in a bounded set. By differentiation we also have that
$$
\pp_q I (q) = \pp_q J_\ve (W_q)  +  \int_{\R^N} \pp_q\Phi(-\Delta)^s\Phi  + V\Phi\pp_q\Phi\ +
$$
$$
\int_{\R^N}
[(W_q +\Phi)^{p} - W_q^{p} - p W_q^{p-1}\Phi] \pp_q W_q  + [(W_q +\Phi)^{p} -   W_q^{p}]\pp_q\Phi .
$$
Since we also have $\|\pp_q\Phi\|_Y= O(\tau)$, then the second and third term above are of size $O(\ve^2)$. Thus,
$$
\pp_q I (q) = \pp_q J_\ve (W_q) +  O(\rho^2) .
$$
uniformly on $q\in \Gamma$
and the proof is complete. \qed

\bigskip

Next we estimate  $ J_\ve (W_q)$ and
$\pp_q J_\ve (W_q)$. We begin with the simpler case $k=1$. Here it is always the case that $$\|E\|_Y + \|\pp_q E\|_Y \le \tau. $$
Let us also set $\xi = \ve q$
We have now that

$$
W_q(x) =  w_\la (x-q) , \quad \la = V(\xi).
$$
We compute
$$
J_\ve (W_q) =   J^\la (w_\la)  + \frac 12 \int_{\R^N} (V(\xi + \ve x)- V(\xi) ) w_\la^2(x)\, dx
$$
where
$$
J^\la (v) =
\frac 12 \int_{\R^{N}} v (-\Delta)^s v   + \frac \la 2  \int_{\R^N} v^2 - \frac 1{p+1}\int_{\R^N} v^{p+1}.
$$

\bigskip
We recall  that
\begin{equation*}
w_\la (x):= \la^{\frac 1{p-1}} w ( \la^{\frac 1{2s}}x)
\end{equation*}
satisfies the equation
$$
(-\Delta)^sw_\la  + \la w_\la - w_\la^p = 0\inn \R^N.
$$
where $w=w_1$  is the unique radial least energy solution of
$$
(-\Delta)^sw  + w - w^p = 0\inn \R^N.
$$
Then, after a change of variables we find
$$
J^\la(w_\la) =
\frac 12 \int_{\R^{N}} w_\la (-\Delta)^s w_\la   + \frac \la 2  \int_{\R^N} w_\la^2 - \frac 1{p+1}\int_{\R^N} w_\la^{p+1} =
 \la^{\frac {p+1}{p-1} -\frac N{2s}} J^1(w).
$$
Now since $w$ is radial, we find
$$ \int_{\R^N} x_i w_\la(x)dx  =0. $$
Thus,
$$
\int_{\R^N} (V(\xi + \ve x)- V(\xi) ) w_\la^2(x)\, dx  = \nn V(\xi)\cdot \int_{\R^N}  xw_\la   + O( \ve^2)   = O(\ve^2)
$$
On the other hand
$$
\pp_q  \int_{\R^N} (V(\xi + \ve x)- V(\xi) ) w_\la^2(x)\, dx =
$$
$$
\ve\int_{\R^N} ( \nabla V(\xi + \ve x)- \nn V(\xi) ) w_\la^2(x)\, dx\  + $$
$$
 2\int_{\R^N} (V(\xi + \ve x)- V(\xi) ) w_\la \pp_q w_\la \, dx \ =  \ O(\ve^2) .
 $$

\begin{lemma}\label{1peak}
Let $\theta =  \frac {p+1}{p-1} -\frac N{2s}$, $c_* = J_1(w)$ and $k=1$.
Then the following expansions hold:
$$
I(q)  =  c_*V^\theta (\xi )  + O(\ve^{\min (4s, 2)})
$$
$$
\nn_q I(q) =  c_*\ve \nabla_\xi  (V^\theta) (\xi)    + O(\ve^{\min (4s, 2)} )  .
$$
\end{lemma}
For the case $k>1$ and $\min_{i\ne j}|q_i -q_j| \ge R>>1$, we observe that, also, $\|E\|_Y = O(\tau)$ and hence
we also have
$$
I (q) = J_\ve (W_q) +  O(\tau^2),  \quad \pp_q I (q) = \pp_q J_\ve (W_q) +  O(\tau^2)  .
$$
By expanding $I(q)$ we get the validity of the following estimate.
\begin{lemma}
\label{kpeak}
Letting $\xi = \ve q$ we have  that
$$
I(q )  =  c_* \sum_{i=1}^k V^\theta (\xi_i ) - \sum_{i \not = j} \frac{c_{ij}}{|q_i -q_j|^{N+2s}} + O(\ve^{\min (4s, 2)} + \frac{1}{R^{ 2(N+2s-\mu)}}),
$$
$$
\nn_q I (q) =  c_* \ve  \nn_\xi\, \big[  \sum_{i=1}^k  V^\theta (\xi_i ) -\sum_{ i \not = j} \frac{c_{ij}}{|q_i- q_j|^{N+2s}}\big ] + O(\ve^{\min (4s, 2)} + \frac{1}{R^{ 2(N+2s-\mu)}})
$$
where $ c_{*}$ and $c_{ij}=  c_0 (V(\xi_i))^{\alpha} (V(\xi_j))^{\beta}$ are positive constants.
\end{lemma}

 \proof
 It suffices to expand $J_\ve(W_q)$. We see that, denoting $w_i(x): = w_{\la_i }(x- q_i)$
 \begin{eqnarray}
J_\ve(W_q) &=&  J_\ve\Big(
 \sum_{i=1}^k  w_{i} \Big)  \ =\  \sum_{i=1}^k J_\ve (w_i)\  +
 \nonumber \\
 & &  \frac 12 \sum_{i \ne j} \int_{\R^N} w_i(-\Delta)^s w_j  + \int_{\R^N}
 V(\ve x) w_i w_j   \nonumber \\
& &  - \frac 1{p+1} \int_{\R^N} \Big(
 \sum_{i=1}^k  w_{i} \Big)^{p+1} -   \sum_{i=1}^k  w_{i}^{p+1}.
 \label{est6}
\end{eqnarray}

We estimate for $ i \ne j$
\begin{eqnarray}
  & & \int_{\R^N} w_i(-\Delta)^s w_j  + \int_{\R^N}
 V(\ve x) w_i w_j  \nonumber \\
 & = & \int_{\R^N} w_i w_j^p + \int_{\R^N}
 (V(\ve x)-\la_j) w_i w_j  \nonumber   \\
&=& (c_{ij}+o(1)) \frac{1}{|q_i-q_j|^{N+2s}} + O(\frac{\ve^{2s}}{R^{N+2s-\mu}})
\label{est7}
\end{eqnarray}
where $c_{ij}= c_0(V(\xi_i))^{\alpha} ((V(\xi_j))^{\beta})$ and $c_0,\alpha, \beta$ are constants depending on $p,s$ and $N$ only.
Indeed,
$$  w_i(x)= \la_i^{\frac{1}{p-1}} w(\la_i^{\frac{1}{2s}} (x-q_i))$$
and it is known that $$ w(x) =  \frac{c_0}{|x|^{N+2s}} (1+ o(1)) \ass |x|\to \infty .
$$
Then, we have
  $$\int_{\R^N} w_j^p w_i
\ =\ \la_{i}^{\frac{1}{p-1}-\frac{n+2s}{2s}}
\la_j^{\frac{p}{p-1}-\frac{n}{2s}} \left(\int_{\R^N} w^p \right )\frac{c_0}{|q_i-q_j|^{N+2s}}, $$
and hence
$$ c_{ij}= c_0 \la_i^{\alpha} \la_j^\beta$$
where  $$ \la_i= V(\xi_i),\quad \la_j=V(\xi_j) , \quad
\alpha= \frac{1}{p-1}-\frac{n+2s}{2s}, \quad\beta= \frac{p}{p-1}-\frac{n}{2s}.
$$

To estimate the last term we note that
\begin{eqnarray}
 \int_{\R^N} \Big( (\sum_{i=1}^k  w_{i} )^{p+1} -   \sum_{i=1}^k  w_{i}^{p+1} \Big)^{p+1} &=& \sum_{j=1}^k \int_{\Omega_j} \Big( (\sum_{i=1}^k  w_{i} )^{p+1} -   \sum_{i=1}^k  w_{i}^{p+1} \Big)^{p+1}
\nonumber \\
&= & \sum_{j=1}^k \sum_{\Omega_j} \Big( (p+1) w_j^p (\sum_{ i \ne j} w_i) + O( w_{j}^{\min (p-1, 1)}  (\sum_{ i \ne j} w_i)^2) \Big)
\nonumber\\
&=& \sum_{j=1}^K \sum_{i \ne j} (p+1) \int_{\R^N} w_j^p w_i + O( \frac{1}{R^{2 (N+2s-\mu)}})
\nonumber \\
&=& \sum_{j=1}^K \sum_{i \ne j} (p+1) \frac{c_{ij}+o(1)}{|q_i-q_j|^{N+2s}} + O( \frac{1}{R^{2 (N+2s-\mu)}})
\label{est8}
\end{eqnarray}

 Substituting (\ref{est7}) and (\ref{est8}) into (\ref{est6}) and using the estimate of $J_\ve (w_i)$ in the proof of Lemma \ref{1peak}, we have estimated $J_\ve (w_i)$, and we have proven the lemma.

\section{The Proofs of Theorems 1--3}
\setcounter{equation}{0}

Based on the asymptotic expansions in Lemma \ref{kpeak}, we  present the proofs of Theorems \ref{teo1}-\ref{teo3}.

{\bf Proof of Theorems \ref{teo1} and \ref{teo2}.}
Let us consider the situation in Remark \ref{remark}, which is more general than that of Theorem 1.
Then, in the definition of the configuration space $\Gamma$ \equ{const},
 we can take a fixed $\delta$ and  $R\sim \ve^{-1}$ and achieve that $ \Lambda \subset \ve\Gamma$. Then
 we get $$\|E\|_Y + \|\pp_q E\|_Y = O(\ve^{\min\{2s,1\}}).$$
Letting
$$ \ttt I(\xi) := I (\ve q)  $$
we need to find a critical point of $\ttt I$ inside $\Lambda$.
By Lemma \ref{kpeak}, we see
then that

$$
\ttt I(\xi) - c_*\vp(\xi) = o(1),   \quad \nn_\xi \ttt I(\xi) - c_*\nn_\xi \vp(\xi) = o(1),
$$
uniformly in $\xi\in \La$ as $\ve\to 0$,
where $\vp$ is the functional in \equ{vp}. It follows, by the assumption on $\vp$  that for all $\ve$ sufficiently small
there exists a $\xi^\ve\in\La$ such that  $\nn \ttt I(\xi^\ve) =0$, hence Lemma \ref{reduction} applies and the desired result
follows.

Theorem \ref{teo2} follows in the same way. We just observe that because of the $C^1$-proximity, the same variational characterization
of the numbers $c$, for the functional $\ttt I(\xi)$ holds. This means that the critical value predicted in that form is indeed close to $c$.
The proof is complete. \qed

\bigskip
{\bf Proof of Theorem 3.}
Finally we prove Theorem \ref{teo3}. Following the argument in \cite{kw}, we choose the following configuration space
\begin{equation}
\Lambda=\{ (\xi_1,..., \xi_k)\ /\ \xi_j \in \Gamma, \min_{i \not = j} |\xi_i-\xi_j| >\varepsilon^{1-\frac{s}{4} } \}
\end{equation}
with $\Gamma$ given by \equ{const},
 and we prove the following Claim and then  Theorem \ref{teo3} follows from Lemma \ref{reduction}:

\medskip

\noindent
{\bf Claim:}  letting $\xi=\ve q$, the problem
\begin{equation}
\max_{ (\xi_1, ..., \xi_k) \in \Lambda} I(q)
\end{equation}
admits a maximizer $(\xi_1^\ve, ..., \xi_k^\ve) \in {\Lambda}.$

\medskip
We shall prove this by contradiction.
First, by continuity of $I(q)$, there is a maximizer  $\xi^\ve =(\xi_1^\ve, ..., \xi_k^\ve) \in \bar{\Lambda}$. We need to prove that
 $\xi \in \La$.  Let us  suppose, by contradiction,
  that $ \xi^\ve \not\in \Lambda$, hence it lies on its boundary. Thus there are two possibilities:
either there is an index $i$ such that $\xi_i^k \in \partial \Gamma$, or there exist indices $ i \ne j$ such that
$$ |\xi_i^\ve-\xi_j^\ve|= \min_{i \ne j} | \xi_i -\xi_j|= \epsilon^{1-s}.$$
Denoting $ q^\ve= \frac{\xi^{\ve}}\ve$, and using Lemma \ref{kpeak}, we have in the first case that
$$
I(q^\ve) \leq c_{*} V^\theta (\xi_i^\ve) + c_{*}\sum_{j \ne i} V^\theta (\xi_j^\ve) + C \ve^{2s}\  \leq
$$
\be\label{est21}
c_{*} k \max_{\Gamma} V^\theta (x) + c_{*} (\max_{\partial \Gamma} V^\theta (x) - \max_{\Gamma} V^\theta (x) )+ C\ve^{2s}.
\ee
In the second case, we invoke again Lemma \ref{kpeak} and obtain
\begin{equation*}
I(q^\ve) \leq c_{*}   k \max_{\Gamma} V^\theta (x)  - c_2 \ve^{\frac{s}{4}} + C\ve^{2s}
\end{equation*}
for some $c_2 >0$.
On the other hand, we can get an upper bound for $ I(q^\ve)$ as follows. Let us choose a point $ \xi_0$ such that
$ V(\xi_0)= \max_{ \Gamma } V(x)$
 and let $$ \xi_j= \xi_0 + \varepsilon^{1- \frac{1}{8} s} (1, 0, ..., 0),\quad  j=1,..., k.$$
  It is easy to see that $(\xi_1, ..., \xi_k) \in \Lambda$. Now, we compute by Lemma \ref{kpeak}:
\begin{equation}
\label{est23}
I(q^\ve)=\max_{\Lambda} I(q) \geq c_{*} k \max_{\Gamma} V^\theta (x) - c_3 \ve^{\frac{s}{8} }.
\end{equation}
For $\ve$ sufficiently small, a contradiction follows immediately from (\ref{est21})-(\ref{est23}).

\qed

\medskip

\noindent
{\bf Acknowledgments}

$\bullet$ J.D. and M.D. have been supported by
 Fondecyt grants  1130360, 110181 and Fondo Basal CMM. J.W. was supported by Croucher-CAS Joint Laboratory and NSERC of Canada.

\medskip

$\bullet$ After completion of this work we have learned about the paper \cite{xxx} in which
 the result of Corollary \ref{coro1} is found for $k=1$ under
further constraints in the space dimension $N$ and the values of $s$ and $p$.

\end{document}